\documentclass[submission]{eptcs}
\providecommand{\event}{QAPL 2011} % Name of the event you are submitting to
\usepackage{breakurl}             % Not needed if you use pdflatex only.

\usepackage{latexsym,amssymb,amsmath,theorem,stmaryrd}

\usepackage{balance}
\usepackage{paralist}
\usepackage{graphicx}
\usepackage{multirow}
 \usepackage{times}
\usepackage[usenames]{color}
\usepackage{gastex}
\usepackage{multirow}
\usepackage{multicol}

\usepackage{subfigure}

\newcommand{\mRNA}{\ensuremath{M}}
\newtheorem{definition}{Definition}
\newtheorem{example}{Example}
\renewcommand{\Pr}[1]{\mathit{Pr}\!\left(#1\right)}
\newcommand{\tmax}{t_{\max}}
\newenvironment{algorithm} 
{ 
	 \begin{figure} }
	{ \end{figure}  
}

\title{On-the-fly Uniformization of Time-Inhomogeneous Infinite Markov Population Models}

\author{
	Aleksandr Andreychenko \qquad Pepijn Crouzen \qquad Linar Mikeev \qquad Verena Wolf
	\institute{Saarland University\\ Saarbr\"{u}cken, Germany} 
	\email{ \{aleand, crouzen, mikeev, wolf\}@cs.uni-saarland.de }
}
\def\titlerunning{On-the-fly Uniformization of Time-Inhomogeneous Infinite Markov Population Models}
\def\authorrunning{A. Andreychenko, P. Crouzen, L. Mikeev, V. Wolf}
\begin{document}
 \maketitle

	\begin{abstract}
		This paper presents an on-the-fly uniformization technique for the analysis
		of time-inhomogeneous Markov population models. This technique is
		%applicable to models with infinite state space and unbounded rates, which
		applicable to models with infinite state spaces and unbounded rates, which
		are, for instance, encountered in the realm of biochemical reaction networks.
		To deal with the infinite state space, we dynamically maintain a finite subset
		of the states where most of the probability mass is located. This approach
		yields an under-approximation of the original, infinite system. We present
		experimental results to show the applicability of our technique.
	\end{abstract}

	\section{Introduction}
		% \fixme[margin]{
		% Markov population models (MPMs) (as Kurtz considers them!):
		% }
		Markov population models (MPMs) are continuous-time 
		Markov processes, where the state of
		the system is a vector of natural numbers (i.e., the populations).
		% \fixme[margin]{applications?}
		Such models are used in various application domains: biology, where the state
		variables describe the population sizes of different organisms, queueing
		theory, where we model a state as a vector of queue occupancies, chemistry,
		where the state variables represent the amount of molecules of different
		chemical species, etc~\cite{Kingman}.
		
		% \fixme[margin]{Measures}
		Besides the expectations and variances of the different populations,
		the probabilities of certain events occurring can be of interest when
		studying MPMs.
		It may be necessary to know the probability of the extinction of a species,
		the probability that a population reaches a certain threshold, or even
		the full distribution of the MPM at a certain time-point, for instance to
		calibrate model parameters.
		% are of interest (e.g. prob. of extinction, reaches threshold,...)
		% or even full distribution up to certain accuracy (for parameter
		% calibration)
		
		% \fixme[margin]{infinite}
		Many Markov population models have \emph{infinitely} many states. In the case
		of biological or chemical applications, we normally cannot provide hard
		upper bounds for population numbers and in the field of queueing theory it
		may be interesting to consider unbounded queues. The evaluation of
		infinite MPMs through numerical \cite{HIBI09} or statistical \cite{gillespie76} 
		analysis has been well-studied for \emph{time-homogeneous} models where
		the dynamics of the system are independent of time.
In \cite{HIBI09} the state space of the model is generated and truncated
on-the-fly during the transient solution, that is, during a certain time
interval only states that are relevant at that time are considered. Thus,
states are added at a certain step and dropped at a later time when they become
irrelevant. A similar technique is proposed in~\cite{Ciardotrunc} for the solution of 
time-homogeneous discrete-time Markov chains.
Note that this is different from on-the-fly techniques for the
computation of steady-state probabilities where the relevent part of the state
space is generated but states are never dropped as time
progresses~\cite{SilvaO92}.

		% \fixme[margin]{time}
Many  Markov models are \emph{time-inhomogeneous}, that is,
		their dynamics change over time. For instance, when
		modeling an epidemic, we may have to take into account that infection rates
		vary seasonally. For traffic models, time-dependent arrival rates can
		be used to model the morning and evening rush hours. In cellular biology
		we see that reaction propensities depend on the cell volume, which waxes
		and wanes as the cell grows and divides. The class of 
		finite time-inhomogeneous Markov models has also been studied in recent years
		 \cite{buchholz,van1992uniformization,MoorselW98}.
		
		In this paper, we develop a numerical 
		algorithm to approximate transient probability
		distributions (i.e., the probability to be in a certain state at a certain
		time) for infinite time-inhomogeneous MPMs. 
		We consider MPMs with state-dependent rates and do not require the
		existence of an upper-bound for the transition rates in the MPM. 
		% To the best of our knowledge no numerical yet exists.
		% \fixme[inline]{Is this true? Mention Linar-Kutta
		% or Engblom or time-inhomogeneous simulation?}
		
		Our algorithm is based on the \emph{uniformization} technique, which is
		a well-known method to approximate the transient probability distribution of
		finite time-homogeneous Markov models~\cite{jensen,Grassmann}. Recently, two
		adaptations of uniformization have been developed. These adaptations
		respectively approximate the transient probabilities for 
		finite time-inhomogeneous~\cite{buchholz} and infinite
time-homogeneous~\cite{HIBI09} Markov models.
		Our algorithm combines
		and refines these two techniques such that
		 infinite time-inhomogeneous MPMs with unbounded rates
		can be tackled. We present two case studies to investigate the 
		effectiveness of
		our approach.

	\section{Markov Population Models}
		Markov chains with large or even infinite 
		state spaces are usually described by some high-level modeling 
		formalism that allows the generation of a (possibly infinite) 
		set of states and transitions.
		Here, we use transition classes to specify a
		Markov population model, that is, a continuous-time 
		Markov chain (CTMC)
		$\left\{X(t),t\ge 0\right\}$  with state space
		$S=\mathbb Z_+^n=\{0,1,\ldots\}^n$, where the $i$-th state variable
		represents the number of instances of the $i$-th species.  Depending
		on the application area, ``species'' stands for types of system
		components, molecules, customers, etc.  The application areas that we
		have in mind are chemical reaction networks, performance evaluation of
		computer systems, logistics, epidemics, etc~\cite{Kingman}.
 
		\begin{definition}[Transition Class] 
		\label{def:traclass}
			A transition class $\tau$ is a triple $(G,w,\alpha)$ where 
		  	$G\subseteq \mathbb Z_+^n$ is the \emph{guard},  
		 	$w\in\mathbb Z^n$ is the \emph{change vector}, and 
		 	$\alpha:G\times \mathbb R_{\ge 0} \to \mathbb R_{\ge 0}$ is 
			the time-dependent \emph{rate function}. 
			Moreover, for any $x\in \mathbb Z_+^n$, we have that $x\in G$ implies $x+w\in \mathbb Z_+^n$.
		\end{definition}
		The guard is the set of states where an instance of $\tau$
		is possible, and if the current state is $x\in G$ then 
		$x+w\in \mathbb Z_+^n$ is the state after an 
		instance of $\tau$ has occurred.
		The rate $\alpha(x,t)$ determines the time-dependent 
		transition probabilities for an infinitesimal time-step $dt$
		$$\Pr{X(t+dt)=x+w\mid X(t)=x}=\alpha(x,t)\cdot dt + o(dt),$$
		where $o$ is a function such that $o(0)=0$ and 
		$\lim_{h\to 0}o(h)/h=0$.
		
		A CTMC $X$ can be specified by a set of $m$ transition 
		classes $\tau_1,\ldots,\tau_m$ as follows.
		For $j\in\{1,\ldots,m\}$, let $\tau_j=(G_j,w_j,\alpha_j)$. 
		For each $t\in\mathbb R_{\ge 0}$ 
		we define the generator matrix $Q(t)$ of $X$ such that 
		the row that describes the transitions of a 
		state $x$ has entry $\alpha_j(x,t)$ at position 
		$Q(t)_{x,x+w_j}$ whenever $x\in G_j$ and zero otherwise.
		Moreover, the diagonal entries of $Q(t)$ are the negative
		sums of the off-diagonal row entries because the    
		row sums of a generator matrix are zero. 
		We assume that each change vector $w_j$ has at least one 
		non-zero entry. To simplify the presentation we assume that 
		all change vectors are distinct.
		We remark that $X$ is called \emph{time-homogeneous} 
		when $Q(t)$ is equal for all $t$. Otherwise, $X$ 
		is called \emph{time-inhomogeneous}. 
		\begin{example}
		\label{ex:gene}
			We consider a simple gene expression model for 
			E. coli cells~\cite{Thattai}. It consists of the transcription of a gene 
			into messenger RNA (mRNA) and subsequent translation of the latter 
			into proteins. 
			%An illustration is given in Fig.~\ref{fig:gene}.
			A state of the 	system is uniquely determined by the
			number of mRNA and protein molecules, that is,
			%a state is a pair $(x_R,x_P)\in\mathbb Z_+^2$.
			a state is a pair $(x_\mRNA,x_P)\in\mathbb Z_+^2$.
			We assume that initially there are no mRNA molecules and no proteins
			in the system, i.e., $\Pr{X(0)=(0,0)}=1$.
		 	Four types of reactions occur in the system.
		 	Let $j\in\{1,\ldots,4\}$ and $\tau_j=\left(G_j,w_j,\alpha_j\right)$
		 	be the transition class that describes the $j$-th reaction type. 	
		 	We first define the guard sets $G_1,\ldots,G_4$ and the change
		 	vectors $w_1,\ldots,w_4$. 
		 	\begin{itemize}
				 \item Transition class $\tau_1$ models gene
					transcription. The corresponding 
					stoichiometric equation is $\emptyset \to$ mRNA.
			 		If a $\tau_1$-transition occurs,
					the number of mRNA molecules increases by one.
					Thus, $w_1=(1,0)$.
					%$u_1(x_R,x_P)=(x_R+1,x_P)$.
					This transition class is possible in all states, i.e.,
					$G_1=\mathbb Z_+^2$. 
				\item We represent the translation of mRNA into 
				 	protein by $\tau_2$ (mRNA $\to$ mRNA+P).
					A $\tau_2$-transition is only possible if
					there is at least  one mRNA molecule in the system. 
				 	We set $G_2=\{(x_\mRNA,x_P)\in\mathbb Z_+^2\mid x_R>0\}$ and
				 	$w_2=(0,1)$.
					%$u_2(x_R,x_P)=(x_R,x_P+1)$. 
					Note that in this case mRNA is
					a reactant that is not consumed.
				 \item Both mRNA and protein molecules can degrade, which is
				 	modelled by $\tau_3$ and $\tau_4$ (mRNA~$\to \emptyset$
					and P~$\to \emptyset$).
					Hence, $G_3=G_2$, $G_4=\{(x_\mRNA,x_P) 
					\in\mathbb Z_+^2\mid x_P>0\}$, 
					$w_3=(-1,0)$, and $w_4=(0,-1)$.
					%$u_3(x_R,x_P)=(x_R-1,x_P)$, and $u_4(x_R,x_P)=(x_R,x_P-1)$.					
			\end{itemize}
		 
			Let $k_1, k_2, k_3, k_4$ be real-valued positive constants.
			We assume that 	transcription  happens
			at rate $\alpha_1(x_\mRNA,x_P,t)= {k_1}\cdot{V(t)}$, 
			that is, the rate is proportional to the cell volume $V(t)$~\cite{wolkenhauer04}.
		 	The (time-independent) translation rate depends 
		 	linearly on the number of mRNA molecules.
			Therefore, $\alpha_2(x_\mRNA,x_P,t)$ $=k_2 \cdot x_\mRNA$.
			%$\alpha_2(x_R,x_P)=\frac{\ln 2}{60} \cdot x_R$. 
		 	%$\alpha_3(x_R,x_P)=\frac{\ln 2}{120}  \cdot x_R$, 
		 	Finally, for degradation, we set
		 	$\alpha_3(x_\mRNA,x_P,t)=k_3\cdot x_\mRNA$ and 
			% $\alpha_4(\vec x)=\frac{\ln 2}{3600} \cdot x_P$,
			$\alpha_4(x_\mRNA,x_P,t)=k_4 \cdot x_P$.
		\end{example}

		We now discuss the transient probability distribution of a MPM.
		Let $S$ be the state space of $X$ and let 
		the transition function $P{(t,t+\Delta)}$ be such that
		the entry for the pair $(x,y)$ of states equals
		$$
		P{(t,t+\Delta)}_{xy}=\Pr{X(t+\Delta)=y\mid X(t)=x}, \quad  t,\Delta\ge 0.
		$$
		If the initial probabilities $\Pr{X(0)=x}$ are specified for each $x\in S$,
		the transient state probabilities $p^{(t)}(x):=\Pr{X(t)=x}$,
		are given by
		$$
		p^{(t)}(y)=\sum\nolimits_{x\in S}p^{(0)}(x)\cdot P{(0,t)}_{xy}.
		$$ 
		We assume that a transition class description 
		uniquely specifies a CTMC and rule out ``pathological cases'' by
		assuming that the sample paths
		%$X{(t)}(\omega)$ are right-continuous step functions. 
		$X{(t)}$ are right-continuous step functions. 
		In this case the transition functions 
		are the unique solution of the
		Kolmogorov backward and forward  equations
		\begin{align}
		&\frac{d}{dt}P(t_0,t)=  Q(t)\cdot P(t_0,t)\\
		&\frac{d}{dt}P(t_0,t)= P(t_0,t)\cdot Q(t), \label{eq:Kolmo}
		\end{align}
		where $0\le  t_0\le t$.
		Multiplication of Eq.~\eqref{eq:Kolmo} with the 
		row vector $p^{(t_0)}$ with entries $p^{(t_0)}(x)$ gives 
		\begin{equation}\label{eq:ODE}
		\frac{d}{dt} p^{(t)}= p^{(t)}\cdot Q(t).
		\end{equation}
		If $S$ is finite, algorithms for the computation of  $p^{(t)}$ 
		are usually based on the numerical integration of
		the linear system of differential equations 
		in Eq.~\eqref{eq:ODE}
		with initial condition $p^{(0)}$.
		Here, we focus on another approach called
		uniformization that is widely used for time-homogeneous 
		Markov chains~\cite{jensen}. It has been adapted  for time-inhomogeneous 
		Markov chains by Van Dijk~\cite{van1992uniformization} and 
		subsequently improved~\cite{MoorselW98,buchholz}.  
		The main advantage of solution techniques based on 
		uniformization is that they provide an underapproximation 
		of the vector $p^{(t)}$ and, thus, provide tight  error bounds.
		Moreover, they are numerically stable and often superior to numerical integration 
		methods in terms of running times~\cite{stewart}.

	\section{Uniformization}
	\label{sec:unif}
		Uniformization is based on the idea to construct, for a CTMC $X$, a Poisson
		process ${N(t),t\ge 0}$ and a subordinated discrete-time 
		Markov chain (DTMC) ${Y(i),i\in\mathbb N}$ such that
		for all $x$ and for all $t$ 
		\begin{equation}\label{eq:unif}
			%\Pr{X(t)=x}=\Pr{Y(N(t))=x}.
			\Pr{X(t)=x}=\Pr{Y(N(t))=x}.
		\end{equation}
		Since Poisson process $N$ and DTMC $Y$ are independent, the equation
		above can be written as
		\begin{equation}\label{eq:unif_sum}
			\Pr{Y(N(t))=x}= \sum_{i=0}^\infty \Pr{Y(i)=x} \Pr{N(t)=i}.
		\end{equation}
		For a finite time-homogeneous MPM with state space $S$ the
		rate $\Lambda$ of the Poisson process $N$ (also called the \emph{uniformization
		rate}) is chosen to be greater than or equal to the maximal exit-rate appearing
		in $X$
		$$ \Lambda \geq \max\limits_{x \in S} \sum_{j=1}^m \alpha_j(x). $$
		For the DTMC $Y$ we find transition probabilities
		$$ \Pr{Y(i\!+\!1)=x\!+\!w_j \mid Y(i) = x} = \frac{\alpha_j(x)}{\Lambda}. $$
		When $X$ is time-inhomogeneous, Arns et al.~\cite{buchholz}
		suggest to define the time-dependent uniformization rate $\Lambda(t)$
		%of $N$ as
		of the \emph{inhomogeneous} Poisson process (IPP) $N$ as
		\begin{equation}\label{eq:maxlambda}
			\Lambda(t) \geq \max_{x\in S}\sum_{j=1}^m \alpha_j(x,t).
		\end{equation}
		% In order to define the transition probabilities of $Y$,
		% we split $\alpha_j(x,t)$ into the two factors $\lambda_j(t)$
		% and $r_j(x)$ such that $\alpha_j(x,t)=\lambda_j(t)\cdot r_j(x)$
		% for all $j,x,t$\footnote{Note that this decomposition is not possible in general, but
		% we assume here that the rate functions of the transition classes are
		% such that they allow such a decomposition. In particular, this decomposition
		% is always possible for chemical reaction networks where the time-dependence
		% stems from fluctuations in reaction volume or temperature.}. Thus, the functions 
		% $\lambda_j:\mathbb R_{\ge 0}\to \mathbb R_{\ge 0}$ contain the 
		% time-dependent part (but are state-independent) and the functions
		% $r_j:S\to \mathbb R_{\ge 0}$ contain the state-dependent part (but are
		% time-independent). 
		For the (time-dependent) transition probabilities 
		of the DTMC $Y$ we then have that $\frac{\alpha_j(x,t)}{\Lambda(t)}$
		is the probability to enter state $x+w_j$ from state $x$ if a state-change
		occurs at time $t$. Arns et al. prove that Eq.~\eqref{eq:unif} holds 
		if the $\alpha_j$ are (right or left) continuous functions in $t$
		and if $S$ is finite (see Theorem 7 in~\cite{buchholz}).
		Here, we relax the latter condition and allow $S$ to be infinite.
		If $\sup_{x\in S} \sum_j \alpha_j(x,t)<\infty$ during 
		the time interval of interest, the proof of Eq.~\eqref{eq:unif} 
		may be expected to proceed along similar lines. In our case, however,
		$\sup_{x\in S} \sum_j \alpha_j(x,t)=\infty$ and
		then the Poisson process $N$ is not well-defined as its rate 
		must be infinite according to Eq.~\eqref{eq:maxlambda}. 
		Therefore, the infinite state 
		space has to be truncated in an appropriate way. 

		\subsection{State Space Truncation}
			We consider a time interval $[t,t+\Delta)$ of length $\Delta$, where the
			transient distribution at time $t$, $p^{(t)}$, of the infinite time-inhomogeneous
			MPM $X$ is known. We now wish
			to approximate the transient distribution at time $t+\Delta$, $p^{(t+\Delta)}$.
			We assume that $p^{(t)}$ has finite support $S_{t,0}$.
			Define $\Pr{N(t,t+\Delta)=i}=\Pr{N(t+\Delta)-N(t)=i}$
			as the probability that $N$ performs $i$ steps within 
			$[t,t+\Delta)$.
			For a fixed positive $\epsilon\ll 1$,
			let $R$ and the rate function $\Lambda$ be such that
			%$S_{t,R}$ is the set of states that are reachable, with probability greater or
			%equal than $1-\epsilon$, from the set $S_{t,0}$ in the time-interval $[t,t+\Delta)$
			$S_{t,R}$ is the set of states that are reachable, with probability greater than
			or equal to $1-\epsilon$, from the set $S_{t,0}$ in the time-interval $[t,t+\Delta)$
			within at most $R$ transitions, i.e.
			\begin{equation}
			\label{eq:truncsum}
				\sum_{i=0}^R \Pr{N(t,t+\Delta)=i} \geq 1-\epsilon.
			\end{equation}
			Furthermore, we have that the rate of $N$ at time $t' \in [t,t+\Delta)$
			must satisfy
			\begin{equation}
			\label{eq:biglambda}
				\Lambda(t') \ge \max_{x\in S_{t,R}}\sum_{j=1}^m \alpha_j(x,t').
			\end{equation}
			Note that $\Lambda(t')$ is adaptive and depends on $t'$, $t$,  $\Delta$, $S_{t,0}$,
			and $R$ as opposed to Arns et al. where $\Lambda(t')$ depends only on $t'$, $t$,
			and $\Delta$ because they consider finite state spaces.
			
			Finding appropriate values for $\Delta$ and $R$ is non-trivial as
			$\Lambda(t')$ determines the speed of the Poisson process $N$ and thereby
			influences the value of $R$. On the other hand, $R$ determines the size of
			the set $S_{t,R}$ and thus influences $\Lambda(t')$.  
			We discuss how to find appropriate choices for $\Delta$ and $R$ given
			the set $S_{t,0}$ in Section~\ref{sec:stepsize}.

			Assume that we find $\Delta$ and $R$ with the above 
			mentioned properties and define $\Lambda(t')$ as in Eq.~\eqref{eq:biglambda}.
			Then, for all $x\in S$, we get an $\epsilon$-approximation
			\begin{equation}\label{eq:approx}
			  \Pr{X(t\!+\!\Delta)\!=\!x}\!\geq\! 
			 \sum_{i=0}^R \Pr{Y(i)\!=\! x\wedge N(t,t\!+\!\Delta)\!=\! i},
			\end{equation}
			where $Y$ has initial distribution $p^{(t)}$.
			The probabilities $\Pr{Y(i)=x\wedge N(t,t+\Delta)=i}$
			can now be approximated in the same way as for the finite
			case~\cite{buchholz}.
			 
			From Eq.~\eqref{eq:approx} we   see that it
			is beneficial if $R$ is small, since this means fewer probabilities
			have to be computed in the right-hand side of Eq.~\eqref{eq:approx}.
			Note that the truncation-point $R$ is small when the uniformization rates
			$\Lambda(t')$ are small during $[t,t+\Delta)$ because if $N$ jumps 
			at a slower rate then  $\Pr{N(t,t+\Delta)>i}$ becomes smaller.
			Thus, it is beneficial to choose
			$\Lambda(t')$ as small as possible while still satisfying Eq.~\eqref{eq:biglambda}.

		\subsection{Bounding approach}
			Let $\hat p^{(t+\Delta)}(x)$ denote the right hand side of
			Eq.~\eqref{eq:approx}, i.e., the approximation of the
			transient  probability of state $x$ at time $t+\Delta$. We compute this
			approximation with the uniformization method as follows.
			The processes $Y$ and $N$
			are independent which implies that  
			%$$\begin{array}{c@{\ }l}
			% &  \Pr{Y(i)\!=\! x\wedge N(t,t\!+\!\Delta)\!=\! i}\\ =&
			%\Pr{Y(i)\!=\! x}\cdot \Pr{N(t,t\!+\!\Delta)\!=\! i}.
			%  \end{array}
			%$$
			$$ \Pr{Y(i)\!=\! x\wedge N(t,t\!+\!\Delta)\!=\! i} = \Pr{Y(i)\!=\! x}\cdot \Pr{N(t,t\!+\!\Delta)\!=\! i}.$$
			The probabilities
			$\Pr{N(t,t+\Delta)=i}$ follow a Poisson distribution with parameter 
			$\bar\Lambda(t,t+\Delta)\cdot \Delta$, where 
			$$\textstyle
			\bar\Lambda(t,t+\Delta)=\frac{1}{\Delta}\int_{t}^{t+\Delta} \Lambda(t')\, dt'.$$
			For the distribution  $\Pr{Y(i)\!=\! x}$, Arns et al. suggest 
			an underapproximation that relies on the fact that for any time-point
			$t' \in [t,t+\Delta)$ we have:
			$$\textstyle
			\frac{\alpha_j(x,t')}{\Lambda(t')}\ge \min_{t''\in[t,t+\Delta]} 
			\frac{\alpha_j(x,t'')}{\Lambda(t'')}=:u_j(x,t,t+\Delta).
			$$
			Thus, for $i\in\{1,2,\ldots,R\}$,
			we iteratively approximate $\Pr{Y(i)\!=\! y}$ as
			\begin{equation}\label{eq:minbound}
			% \begin{array}{l@{\ }c@{\ }l}
			%\hspace{-3ex}
			%\Pr{Y(i)\!=\! y}&\ge&\hspace{-2ex}\sum\limits_{x,j:y=x+v_j} \hspace{-2ex}
			%\Pr{Y(i-1)\!=\! x}\cdot u_j(x,t,t\! +\! \Delta)\\[2ex]
			%&+& \Pr{Y(i-1)\!=\! y}\cdot u_0(y,t,t\! +\! \Delta).
			%\end{array}
			\Pr{Y(i)\!=\! y} \ge \sum\limits_{x,j:y=x+w_j} \Pr{Y(i-1)\!=\! x}\cdot u_j(x,t,t\! +\! \Delta) 
									+ \Pr{Y(i-1)\!=\! y}\cdot u_0(y,t,t\! +\! \Delta).
			\end{equation}
			Here, $x$ ranges over all direct predecessors of $y$ 
			and  the self-loop probability $u_0(y,t,t+\Delta)$
			of   $y$ is given by
			$$\textstyle u_0(y,t,t+\Delta)= \min\limits_{t'\in[t,t+\Delta]} 
			\left(1-\sum\limits_{j=1}^m
			\frac{\alpha_j(y,t')}{\Lambda(t')}\right).$$
			Note that often we   can   split   $\alpha_j(x,t')$ into 
			two factors $\lambda_j(t')$
			and $r_j(x)$ such that $\alpha_j(x,t')=\lambda_j(t')\cdot r_j(x)$
			for all $t',j,x$\footnote{Note that this decomposition
			is always possible for chemical reaction networks where the time-dependence
			stems from fluctuations in reaction volume or temperature.}. 
			Thus, the functions 
			$\lambda_j:\mathbb R_{\ge 0}\to \mathbb R_{> 0}$ contain the 
			time-dependent part (but are state-independent) and the functions
			$r_j:S\to \mathbb R_{> 0}$ contain the state-dependent part (but are
			time-independent). Then each minimum   defined above can be computed 
			for all states by considering 
			$$\textstyle
			 \min\limits_{t'\in[t,t+\Delta]} \frac{\lambda_j(t')}{\Lambda(t')}.
			$$
			In particular, if $\lambda_j$ and $\Lambda$ are monotone, 
			the above minimum is easily found analytically.
			
			The approximation in Eq.~\eqref{eq:minbound} implies that
			for the time interval 
			$[t,t+\Delta)$, we compute a sequence of substochastic vectors
			$v^{(1)},v^{(2)}, \ldots, v^{(R)}$ to approximate 
			%the probabilities $\Pr{Y(i)=x}$. Initially we start the DTMC $Y$ with 
			the probabilities $\Pr{Y(i)=y}$. Initially we start the DTMC $Y$ with 
			the  approximation $\hat p^{(t)}=:v^{(0)}$  
			of the previous step. Then we compute $v^{(i+1)}$ from 
			$v^{(i)}$ based on the transition probabilities 
			$u_j(x,t,t+\Delta)$ for $i\in\{0,1,\ldots, R\}$. 
			Since these transition probabilities may sum up to less than one,
			the resulting vector $v^{(i+1)}$ may also sum up to less than one.
			Since, for the computation of  $\hat p^{t+\Delta}$, we weight these vectors
			with the Poisson probabilities and add them  up the underapproximation $\hat p^{t+\Delta}$
			contains an additional approximation error.
			  In general, the larger
			the time-period $\Delta$, the worse the underapproximations $u_j(x,t,t+\Delta)$
			are and thus the underapproximation $\hat p^{t+\Delta}$ becomes worse as well. 
			We illustrate this effect by applying the bounding approach to our
			running example.
			% \fixme{it is unclear to me how
			% the minima can be found efficiently since the function in Alg.1
			% gives only the average rate of the Poisson process! Or should  
			% the function FindMaxInt give as an output also the state 
			% were the maximum was found? Is it clear that at an earlier time instant
			% no other state has a higher exit rate???}
			\begin{example}
			\label{ex:2}
				In the gene expression of Example~\ref{ex:gene}, the time-dependence 
				is due to the volume and only affects the rate function 
				$\alpha_1$ of the first transition class.
				The time until an E. coli cell divides varies widely from about
				20 minutes to many hours and depends on growth conditions.
				Here, we assume a cell cycle time of one hour and a linear 
				growth~\cite{arkin2}.
				Thus, if at time $t=0$ we consider a cell immediately after 
				division  then the cell volume doubles after   3600 sec. 
				Assume that $\Delta\le 3600$. Then,    
				$\alpha_1(x,t')=k_1'\cdot (1+\frac{t'}{3600})$ for all $x\in S$.
				Assume we have a right truncation point $R$ such that
				%\begin{center}
				$$\Lambda(t')=\max\limits_{x_R,x_P}k_1'\cdot (1+\frac{t'}{3600})+(k_2+k_3)\cdot x_R+k_4\cdot x_P$$
				%\end{center}
				where $x_R$ and $x_P$ range over all states $(x_R,x_P) \in S_{0,R}$ and
				Eq.~\eqref{eq:truncsum} holds.
				Then we find, for each 
				time-point $t' \in [0,\Delta)$, the same state   for which
				the exit-rate $\alpha_0(x,t'):=\sum_{j=1}^m \alpha_j(x,t')$
				is maximal, since the only time-dependent propensity is
				independent of the state-variables. Let  $(x^{\max}_R,x^{\max}_P)$ denote
				this state. In general this is not the case, for instance in the realm
				of chemical reaction systems we have that the propensities of bimolecular
				reactions (reactions of the from $A + B \rightarrow \ldots$) are dependent
				both on cell-volume and the population numbers. For such a system we
				may find that different states have the maximal exit-rate within the
				time-frame $[0,\Delta)$. We discuss how to overcome this difficulty
				in Subsection~\ref{sec:maxrates}.
				The transition probabilities of the DTMC $Y$ are now defined as
				%\begin{align*}
				 $$ u_1(x_R,x_P,0,\Delta) =  \min\limits_{t'\in[0,\Delta]} \frac{\alpha_1(x_R,x_P,t')}{\Lambda(t')}
										=  \frac{\alpha_1(x,0)}{\Lambda(0)}
										=  \frac{k_1'}{k_1' +(k_2+k_3)\cdot x^{\max}_R+k_4\cdot x^{\max}_P}$$
				%\end{align*}
				and, for  $j\in\{2,3\}$,
				%\begin{align*}
				 $$u_j(x_R,x_P,0,\Delta) = \min\limits_{t'\in[0,\Delta]} \frac{\alpha_j(x_R,x_P,t')}{\Lambda(t')}
										= \min\limits_{t'\in[0,\Delta]} \frac{k_j \cdot x_R}{\Lambda(\Delta)}
										=\frac{k_j\cdot x_R}{k_1' \cdot (1+\frac{\Delta}{3600})+(k_2+k_3)\cdot x^{\max}_R+k_4\cdot x^{\max}_P},
				$$
				$$u_4(x_R,x_P,0,\Delta)	=\frac{k_4\cdot x_P}{k_1' \cdot (1+\frac{\Delta}{3600})+(k_2+k_3)\cdot x^{\max}_R+k_4\cdot x^{\max}_P}.$$
				%\end{align*}
				For the self-loop probability we find:
				\begin{align*}
				u_0(x_R,x_P,0,\Delta) &= \min_{t' \in [0,\Delta]}  \left(1 - \sum_{j=1}^4 \frac{\alpha_j(x_R,x_P,t')}{\Lambda(t')} \right)
										= \left( 1 - \max_{t' \in [0,\Delta)}  \sum_{j=1}^4 \frac{\alpha_j(x_R,x_P,t')}{\Lambda(t')} \right) \\
				&\hspace{-45pt}= 1 - \sum_{j=1}^4 \frac{\alpha_j(x_R,x_P,\Delta)}{\Lambda(\Delta)}
				= 1 - \frac{k_1' \cdot (1+ \frac{\Delta}{3600}) + (k_2+k_3) \cdot x_R + k_4 \cdot x_P}
				{k_1' \cdot (1+ \frac{\Delta}{3600}) + (k_2+k_3) \cdot x^{\max}_R + k_4 \cdot x^{\max}_P}.
				\end{align*}
			
				We now calculate the fraction of probability lost during the computation of $v^{(i+1)}$ from $v^{(i)}$,	i.e.,
				\begin{align*}
				&1-\sum_{j=0}^4 u_j(x_R,x_P,0,\Delta) = 
								\frac{k_1' \cdot (1 + \frac{\Delta}{3600})}	
								{k_1' \cdot (1+ \frac{\Delta}{3600}) + (k_2+k_3) \cdot x^{\max}_R + k_4 \cdot x^{\max}_P}
								- \frac{k_1'}{k_1' +(k_2+k_3)\cdot x^{\max}_R+k_4\cdot x^{\max}_P} \\
				&\hspace{45pt} = \frac{(k_2+k_3)\cdot x^{\max}_R+k_4\cdot x^{\max}_P}{k_1' +(k_2+k_3)\cdot x^{\max}_R+k_4\cdot x^{\max}_P}
				- \frac{(k_2+k_3)\cdot x^{\max}_R+k_4\cdot x^{\max}_P}{k_1' \cdot (1+ \frac{\Delta}{3600}) + (k_2+k_3) \cdot x^{\max}_R + k_4 \cdot x^{\max}_P}.
				\end{align*}
				For $\Delta=0$ we have a probability loss of $0$ and for $\Delta>0$ we can
				see that the probability loss increases with increasing $\Delta$.
			\end{example}

		\subsection{Time-stepping approach} 
		\label{sec:timestep}
\begin{figure*}[t]
				\centerline{ 
				\scalebox{0.85}{
				\unitlength=4pt
				\begin{picture}(120,40)(0,0)
					\gasset{AHnb=0}
					\put(0,0){
					  \put(8,32){Support at time $t$}
					  \drawline(5,30)(5,5)(35,5)
					  \put(2,25){$x_2$}
					  \drawline[AHnb=1](3,27)(3,30)
					  \put(29,2){$x_1$}
					  \drawline[AHnb=1](32,3)(35,3)
					  \drawccurve(16,15)(19,17.5)(21,17.5)(24,15)(21,12.5)(19,12.5)
					  \put(18.4,14){$S_{t,0}$}
					}
					\put(40,0){
					  \put(6.5,32){Truncation for the first step}
					  \drawline(5,30)(5,5)(35,5)
					  \put(2,25){$x_2$}
					  \drawline[AHnb=1](3,27)(3,30)
					  \put(29,2){$x_1$}
					  \drawline[AHnb=1](32,3)(35,3)
					  \drawccurve(16,15)(19,17.5)(21,17.5)(24,15)(21,12.5)(19,12.5)
					  \put(18.4,14){$S_{t,0}$}
					  \drawccurve(12,15)(17,21)(23,21)(28,15)(23,9)(17,9)
					  \put(18.4,19){$S_{t,R}$}
					} 
					\put(80,0){
					  \put(6.5,32){Truncation for the second step}
					  \drawline(5,30)(5,5)(35,5)
					  \put(2,25){$x_2$}
					  \drawline[AHnb=1](3,27)(3,30)
					  \put(29,2){$x_1$}
					  \drawline[AHnb=1](32,3)(35,3)
					  \drawccurve(16,15)(19,17.5)(21,17.5)(24,15)(21,12.5)(19,12.5)
					  \put(18.4,14){$S_{t,0}$}
					  \drawccurve(12,15)(17,21)(23,21)(28,15)(23,9)(17,9)
					  \put(16.4,18.5){$S_{t+\Delta,0}$}
					  \drawccurve(8,15)(15,24)(25,24)(32,15)(25,6)(15,6)
					  \put(16.4,22.5){$S_{t+\Delta,R}$}
						}
					\end{picture}
					}
				}
				\caption{Illustration of the state space truncation
				approach for the two-dimensional case. Given the distribution 
				$\hat p^{(t)}$ with support $S_{t,0}$, a truncation point $R$
				and a time-step $\Delta$, we compute in the  first step the 
				distribution $\hat p^{(t+\Delta)}$ with support $S_{t,R}=S_{t+\Delta,0}$.
				For the next step we  consider the set $S_{t+\Delta,R}$.
				}\label{fig:trunc}
			\end{figure*}
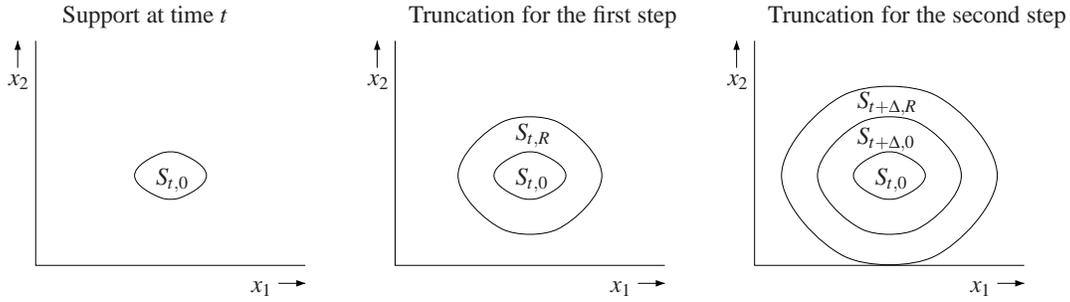

			Given that a large time horizon may lead to decreased accuracy,
			Arns et al.~\cite{buchholz} suggest to partition
			the time period of interest $[0,\tmax)$ in steps of length $\Delta$.
			In each step, an approximation 
			of the transient distribution at the current time instant, $\hat p^{(t)}$,
			is computed and used as initial condition for the next step. 
			The number of states that we consider, that is, $|S_{t,R}|$
			grows in each step. The probabilities of all remaining 
			states of $S$ are approximated as zero.  
			Thus, each step yields a vector $\hat p^{(t+\Delta)}$ with positive 
			entries for all states  $x\in S_{t,R}$ that approximate 
			$\Pr{X(t+\Delta)=x}$. The vector $\hat p^{(t+\Delta)}$ with support
			$S_{t,R}=S_{t+\Delta,0}$ is then used as the initial distribution
			to approximate the vector $\hat p^{(t+\Delta+\Delta')}$. See Figure~\ref{fig:trunc}
			for a sketch of the state truncation approach. Note that the chosen time-period
			$\Delta$ may vary for different steps of the approach.

			It is easy to see that the total error is the sum of  the errors in 
			each step, where the error of a single step equals the amount of probability 
			mass that ``got lost'' due to the underapproximation.
			More precisely, we have two sources of error, namely the
			error due to the truncation of the infinite sum in Eq.~\eqref{eq:unif_sum}
			and the error due to the bounding approach that relies on 
			Eq.~\eqref{eq:minbound}.
			
			  In~\cite{buchholz}, Arns et al. give \emph{exact} formulas for 
			the first three terms of the sum in Eq.~\eqref{eq:approx} (for $i=0,1,2$).
			Thus, if the approximation $\hat p^{(t)}$ of $p^{(t)}$ is exact, then 
			$\hat p^{(t+\Delta)}$ is an underapproximation
			due to the remaining terms  in Eq.~\eqref{eq:approx}.
			This implies that the smaller  $R$ becomes, the closer 
			the error will be to the error bound~$\epsilon$. On the other
			hand, a small truncation point means that only a 
			small time step $\Delta$ is possible (see Eq.~\eqref{eq:truncsum}),
			which means that many steps are necessary until the 
			final time instant $\tmax$ is reached. In order to 
			explore the trade-off between running time and accuracy,
			we run experiments with different values for the 
			predefined truncation point $R$ that determines the 
			step size $\Delta$. We report on these experiments 
			in Section~\ref{sec:expr}.

		\section{On-the-fly Algorithm}\label{sec:onthefly}
			As we can see in Figure~\ref{fig:trunc},
			the number of states that are considered to compute 
			  $\hat p^{(t^{max})}$ from $\hat p^{(t)}$ grows in each step, since
			all states within a radius of $R$ transitions from 
			a state in the previous set $S_{t,0}$ are added. This makes the  
			approach infeasible for Markov models with a large or even 
			infinite state space because the memory requirements 
			are too large. Therefore, we suggest to use a similar 
			strategy as described in previous work~\cite{HIBI09}
			to keep the memory requirements low and achieve faster 
			running times.
			
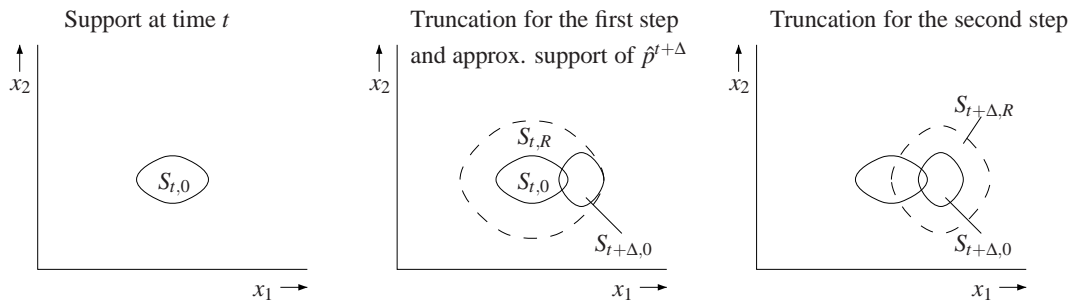
\begin{figure*}[t]
					\centerline{ 
					\scalebox{0.85}{
					\unitlength=4pt
					\begin{picture}(120,40)(0,0)
						\gasset{AHnb=0}
						\put(0,0){
						  \put(8,32){Support at time $t$}
						  \drawline(5,30)(5,5)(35,5)
						  \put(2,25){$x_2$}
						  \drawline[AHnb=1](3,27)(3,30)
						  \put(29,2){$x_1$}
						  \drawline[AHnb=1](32,3)(35,3)
						  \drawccurve(16,15)(19,17.5)(21,17.5)(24,15)(21,12.5)(19,12.5)
						  \put(18.4,14){$S_{t,0}$}
						}
						\put(40,0){
						  \put(6.5,32){Truncation for the first step}
						  \put(6.5,28){and approx. support of $\hat p^{t+\Delta}$}
						  \drawline(5,30)(5,5)(35,5)
						  \put(2,25){$x_2$}
						  \drawline[AHnb=1](3,27)(3,30)
						  \put(29,2){$x_1$}
						  \drawline[AHnb=1](32,3)(35,3)
						  \drawccurve(16,15)(19,17.5)(21,17.5)(24,15)(21,12.5)(19,12.5)
						  \put(18.4,14){$S_{t,0}$}
						  \drawccurve[dash={1.5}0](12,15)(17,21)(23,21)(28,15)(23,9)(17,9)
						  \put(18.4,19){$S_{t,R}$}
						  \drawccurve(23,15)(25.5,18)(28,15)(25.5,12)
						  \put(27,7){$S_{t+\Delta,0}$}
						  \drawline(30,9)(26,13)
						} 
						\put(80,0){
						  \put(6.5,32){Truncation for the second step}
						  \drawline(5,30)(5,5)(35,5)
						  \put(2,25){$x_2$}
						  \drawline[AHnb=1](3,27)(3,30)
						  \put(29,2){$x_1$}
						  \drawline[AHnb=1](32,3)(35,3)
						  \drawccurve(16,15)(19,17.5)(21,17.5)(24,15)(21,12.5)(19,12.5)
						  %\put(18.4,14){$S_{t,0}$}
						  \drawccurve(23,15)(25.5,18)(28,15)(25.5,12)
						  \put(27,7){$S_{t+\Delta,0}$}
						  \drawline(30,9)(26,13)
						  \drawccurve[dash={1.5}0](20,15)(25.5,21)(31,15)(25.5,9)
						  \put(27,22.5){$S_{t+\Delta,R}$}
						  \drawline(30,22)(28,19)
						}
					\end{picture}
					}
					}
					\caption{Illustration of the on-the-fly algorithm for the two-dimensional case.
				Given the distribution 
				$\hat p^{(t)}$ with support $S_{t,0}$, a truncation point $R$
				and a time-step $\Delta$, we compute in the  first step the 
				distribution $\hat p^{(t+\Delta)}$ with approximate support
				$S_{t+\Delta,0} \subset S_{t,R}$.
				For the next step we consider the set $S_{t+\Delta,R}$.
				}\label{fig:otf}
				\end{figure*}

			The underlying principle of this approach 
			is to dynamically maintain a snapshot of the part of the state space
			where most of the transient probability distribution is located.
			We achieve this by adding and removing states in an on-the-fly fashion. 
			The decision which states to add and which states to remove 
			depends on a small probability threshold $\delta>0$.
			After the computation of the vector $v^{(i+1)}$ 
			based on $v^{(i)}$, we set all entries in $v^{(i+1)}$ to zero that have a probability 
			less than $\delta$. This significantly reduces the computational complexity 
			since only parts of the transition probability matrix of $Y$ have to be generated~\cite{HIBI09}
			(for instance, we explore $360000$ states at time instant $t=600$
			for the gene expression system of Example~\ref{ex:gene}	if $\delta=0$ but only $5700$
			states are stored when $\delta=10^{-15}$).
			Let
			$$
			S^{(0)}:=\{x:v^{(0)}(x)>0\}=S_{t,0}
			$$
			and, for $i\in\{1,\ldots,R\}$ let $S^{(i)}$ be the set of states that we consider 
			to compute $v^{(i+1)}$ from $v^{(i)}$.   
			We remark that this also decreases the speed   of the Poisson 
			process $N$ since the sets $S_{t,0}$ and $S_{t,R}$ are smaller and thus 
			 the maximum in Eq.~\eqref{eq:biglambda} is now taken 
			over fewer states. We illustrate this effect in Figure~\ref{fig:otf}.
			 This effect is particularly important if during 
			an interval $[t,t^{\max})$ in certain 
			parts of the state space the dynamics of the system is fast 
			while it is slow in other parts where the latter contain the main part of the probability 
			mass. On the other hand, the threshold $\delta$ 
			introduces another approximation error which may become 
			large if the time horizon of interest is long. 
			Moreover, if $\rho$ is
			a bound for the error introduced by the above strategy of neglecting
			%certain states, we can reserve a portion of $\rho\cdot \frac{\Delta}{\tmax}$
			certain states, we can reserve a portion of the probability loss 
			$\rho\cdot \frac{\Delta}{\tmax}$
			for the interval $[t,t+\Delta)$ and repeat the computation with a smaller 
			threshold $\delta$ if more than the allowed portion of probability was
			neglected.
%  Note that we can easily track how much probability got ``lost'' 
% 			by adding up the probability inflow that was not added to any income-field.

			%\aa{Review: better to reformulate like the following.
			%Since all the introduced approximations are underapproximations
			%we are able to compute the total error of the approximation..}
			The approximation that we suggest above is again an underapproximation
			and since the approximations suggested in the  
			previous sections are also underapproximations, we are still able to compute the total  error of the 
			approximation $\hat p^{(t)}$ of $p^{(t)}$ as
			\begin{equation}
			\label{eq:totalerror}
				%\textstyle1-\sum_{x\in S_{t,R}}\hat p^{(t)}(x).
				1-\sum_{x\in S_{t,R}}\hat p^{(t)}(x).
			\end{equation}
			Clearly, $t'>t$ implies that the error at time $t'$ is higher 
			than the error at time $t$.
			For our experimental 
			results in Section~\ref{sec:expr} we choose $\delta=10^{-15}$
			and report on the total error of the approximation at time $\tmax$.

		\subsection{Determining the step-size}\label{sec:stepsize}
			Given an error bound $\epsilon>0$, a
			time-point $t$, for which the support of $\hat p^{(t)}$
			is $S_{t,0}$, and a time-point $\tmax$ for which we wish to approximate
			the transient probability distribution, we now discuss how to find a
			time-step $\Delta$ such that 
			Eqs.~\eqref{eq:truncsum} and~\eqref{eq:biglambda} hold.
			Recall that the probabilities $\Pr{N(t,t+\Delta)=i}$ follow a Poisson 
			distribution with parameter 
			$\bar\Lambda(t,t+\Delta) \cdot \Delta$, which we denote by
			$\mu_{R,\Delta}$ to emphasize the dependence on $\Delta$
			and the right truncation point $R$. Note that the latter dependence is due to the 
			maximum in Eq.~\eqref{eq:biglambda} that is 
			defined over the set $S_{t,R}$, the set of all states that are reachable
			from a state in $S_{t,0}$ by at most $R$ transitions. We have
			\begin{equation} \label{eq:mu}
			\mu_{R,\Delta} = \int_{t}^{t+\Delta} \Lambda(t')\, dt'.
			\end{equation}

			Here, we propose to first choose a desired right truncation point 
			 $R^*$ and then find  a time-step $\Delta$ such that 
			Eqs.~\eqref{eq:truncsum} and~\eqref{eq:biglambda} hold.
			We perform an iteration where in each step we systematically 
			choose different values for $\Delta$ and compare the associated right 
			truncation point $R$ with $R^*$.
			Since $\mu_{R^*,\Delta}$ is monotone in $\Delta$ this can be done in a binary search 
			%fashion as described in Algorithm~\ref{alg:binarysearch}.
			fashion as described in Algorithm~\ref{alg:binarysearch} .
			We start with the two bounds $\Delta^-=0$ and $\Delta^+=\tmax-t$. The function $\textrm{FindMaxState}(\Delta,R^*)$ finds a state $x^{\max}$ such that for
			all time-points $t' \in [t,t+\Delta)$ we have
			\begin{equation}\label{eq:xmax}
			\sum_{j=1}^m \alpha_j(x^{\max},t') \geq \max\limits_{x' \in S_{t,R^*}}
			\sum_{j=1}^m \alpha_j(x',t').
			\end{equation}

			The choice of $x^{\max}$ also determines the uniformization rate
			$$\Lambda(t') = \sum_{j=1}^m \alpha_j(x^{\max},t').$$
			It immediately follows from Eq.~\eqref{eq:xmax} that Eq.~\eqref{eq:biglambda}
			holds. In Section~\ref{sec:maxrates}, we discuss how to find $\Lambda$ efficiently 
			by selecting a state $x^{\max}$, while avoiding
			that the uniformization rates $\Lambda(t')$ are chosen to be very large.
			
			The function $\textrm{ComputeParameter}(t,t+\Delta,x^{\max})$ now computes
			the integral $\mu_{R^*,\Delta}$ using $x^{\max}$. If possible we compute
			the integral analytically, otherwise we use a numerical integration
			technique. The function $\textrm{FoxGlynn}(\mu,\epsilon)$
			computes the right truncation point of a homogeneous Poisson 
			process with rate $\mu$ for a given error bound $\epsilon$,
			i.e. the value $\hat R$ that is the smallest positive integer 
			such that 
			%\begin{center}
			$$ \sum_{i=0}^{\hat R}\frac{\mu^i}{i!}e^{-\mu}\ge 1-\epsilon. $$
			%\end{center}
			For the refinement of the bounds $\Delta^-$ and $\Delta^+$ in lines
			13--17 we exploit that $R$ is monotone in $\Delta$.

%			\setcounter{figure}{0}
%			\begin{center}
%			\begin{figure}[t]
%			  \includegraphics[width= \linewidth]{algs2.eps}
%			\end{figure}
%			\end{center}

			% =========================================================================
			\setcounter{figure}{0}
			\begin{algorithm}[t] \renewcommand\arraystretch{1.3}
				\centering
				\subfigure[The step size $\Delta$ is determined in a 
				binary-search fashion.\label{alg:binarysearch}]
				{
				\scalebox{0.9}{
				\begin{tabular}{|l|l|} 
				\hline 
				Input   &$R^*$, $t$, $\tmax$, $\epsilon$ \\
				\hline
				Output  & $\Delta$, $x^{\max}$\\%$\mu_{R^*,\Delta}$\\ 
				\hline
				Global  & State space $\hat S$, ...\\
				\hline 
				\multicolumn{2}{|l|}{1 \hspace{0ex} $\Delta^+ := \tmax-t;$ \emph{//upper bound for $\Delta$}}\\[-0.5ex] 
				\multicolumn{2}{|l|}{2 \hspace{0ex} $R := 0;$} \\[-0.5ex]
				\multicolumn{2}{|l|}{3  \hspace{0ex} $x^{\max}:=\textrm{FindMaxState}(\Delta^+,R^*);$} \\[-0.5ex] 
				\multicolumn{2}{|l|}{4  \hspace{0ex} $\mu_{R^*,\Delta^+}:= \textrm{ComputeParameter}(t,t+\Delta^+,x^{\max})$ }  
												\\[-0.5ex] 
				\multicolumn{2}{|l|}{5 \hspace{0ex} $R^+ := \textrm{FoxGlynn}(\mu_{R^*,\Delta^+},\epsilon);$} \\[-0.5ex] 
				\multicolumn{2}{|l|}{6 \hspace{0ex} \textbf{if} $R^+\le R^*$ \textbf{then}}  \\[-0.5ex] 
				\multicolumn{2}{|l|}{7 \hspace{2ex} $\Delta := \Delta^+;$}  \\[-0.5ex] 
				\multicolumn{2}{|l|}{8 \hspace{0ex} \textbf{else} }\\[-0.5ex] 
				\multicolumn{2}{|l|}{9 \hspace{3ex} $R^-:=0; \Delta^-:=0;$ \emph{//lower bound for $\Delta$}}\\[-0.5ex] 
				\multicolumn{2}{|l|}{10 \hspace{2ex} \textbf{while} $R\neq R^*$ }\\[-0.5ex] 
				\multicolumn{2}{|l|}{11\hspace{4ex} $\Delta := \frac{\Delta^+-\Delta^-}{2};$}\\[-0.5ex] 
				\multicolumn{2}{|l|}{12\hspace{4ex} $\mu_{R^*,\Delta}:=\textrm{ComputeParameter}(\Delta,R^*);$}\\[-0.5ex] 
				\multicolumn{2}{|l|}{13\hspace{4ex} $R := \textrm{FoxGlynn}(\mu_{R^*,\Delta},\epsilon);$} \\[-0.5ex] 
				\multicolumn{2}{|l|}{14\hspace{4ex} \textbf{if} $R^-<R^*<R$  } \\[-0.5ex] 
				\multicolumn{2}{|l|}{15\hspace{6ex} $R^+:=R;\Delta^+:=\Delta;$}\\[-0.5ex] 
				\multicolumn{2}{|l|}{16\hspace{4ex} \textbf{elseif} $R<R^*<R^+$  } \\[-0.5ex] 
				\multicolumn{2}{|l|}{17\hspace{6ex} $R^-:=R;\Delta^-:=\Delta;$}\\[-0.5ex] 
				\multicolumn{2}{|l|}{18\hspace{4ex} \textbf{endif}    }\\[-0.5ex] 
				\multicolumn{2}{|l|}{19\hspace{2ex} \textbf{endwhile}    }\\[-0.5ex] 
				\multicolumn{2}{|l|}{20\hspace{0ex} \textbf{endif}    }\\[-0.0ex] 
				\hline
				\end{tabular}
				}
				 }
				 \subfigure[The complete algorithm.\label{alg:complete_algorithm}]
				{
					\scalebox{0.9}
					{
				\begin{tabular}{|l|l|} 
				\hline 
				Input   &$t_0$, $\tmax$, $p^{t_0}$, $\epsilon$, $R^*$ \\
				\hline
				Output  & $p^{t_1}$, $p^{t_2}$, ..., $p^{\tmax}$\\
				\hline
				Global  & State space $\hat S$, ...\\
				\hline 
				\multicolumn{2}{|l|}{1 \hspace{0ex} $t_{cur} := t_0;$ }
													\\[-0.5ex] 
				\multicolumn{2}{|l|}{2 \hspace{0ex} $\left( \Delta, x^{\max} \right) := 
													\textrm{Algorithm 1}(R^*,t_{cur},\tmax,\epsilon);$ }
													\\[-0.5ex] 
				\multicolumn{2}{|l|}{3 \hspace{0ex}	$t_{next}:=t_{cur} + \Delta;$}
													\\[-0.5ex] 
				\multicolumn{2}{|l|}{4 \hspace{0ex} $\mu:= \textrm{ComputeParameter}(t_{cur},t_{next},x^{\max});$ }
													\\[-0.5ex] 
				\multicolumn{2}{|l|}{5 \hspace{0ex} \textbf{while} $t_{cur} \le \tmax$ }
													\\[-0.5ex]
				\multicolumn{2}{|l|}{6 \hspace{2ex} $i:=1$ }
													\\[-0.5ex] 
				\multicolumn{2}{|l|}{7 \hspace{2ex} \textbf{while} $i \le R^*$ }
													\\[-0.5ex]
				\multicolumn{2}{|l|}{8 \hspace{4ex} $\textrm{Compute }v^{(i)}(x);$
													\emph{//DTMC probabilities} } 
													\\[-0.5ex] 
				\multicolumn{2}{|l|}{9 \hspace{4ex} $\textrm{Compute IPP N probabilities};$ }
													\\[-0.5ex] 
				\multicolumn{2}{|l|}{10\hspace{4ex} $\textrm{Accumulate }\hat{p}^{t_{cur}}(x);$
													\emph{//CTMC probabilities} } 
													\\[-0.5ex]
				\multicolumn{2}{|l|}{11\hspace{4ex} $i:=i+1;$}	
													\\[-0.5ex]
				\multicolumn{2}{|l|}{12\hspace{2ex} \textbf{endwhile} }
													\\[-0.5ex] 
				\multicolumn{2}{|l|}{13\hspace{2ex} $t_{cur}:=t_{next};$ }
													\\[-0.5ex]
				\multicolumn{2}{|l|}{14\hspace{2ex} $\left( \Delta, x^{\max} \right) := 
													\textrm{Algorithm 1}(R^*,t_{cur},\tmax,\epsilon);$ }
													\\[-0.5ex]
				\multicolumn{2}{|l|}{15\hspace{2ex}	$t_{next}:=t_{cur} + \Delta;$}
													\\[-0.5ex]
				\multicolumn{2}{|l|}{16\hspace{2ex} $\mu:= \textrm{ComputeParameter}(t_{cur},t_{next},x^{\max});$ }
													\\[-0.5ex]
				\multicolumn{2}{|l|}{17\hspace{0ex} \textbf{endwhile} }
													\\[-0.5ex]
				\multicolumn{2}{|l|}{\hspace{0ex}  } \\[-0.5ex]
				\multicolumn{2}{|l|}{\hspace{0ex}  } \\[-0.5ex]
				\multicolumn{2}{|l|}{\hspace{0ex}  }
													\\[-0.0ex]
				\hline
				\end{tabular} 
				}
				}
				\caption{Algorithms}
				%%\vspace{-3ex} 
			\end{algorithm}
			% =========================================================================

			\setcounter{figure}{1}

		\subsection{Determining the maximal rates}\label{sec:maxrates}
			%The function $\textrm{FindMaxState}(\Delta,R^*)$ in Algorithm~\ref{alg:binarysearch}
			The function $\textrm{FindMaxState}(\Delta,R^*)$ in Algorithm~\ref{alg:binarysearch}
			finds a state $x^{\max}$ such that its exit-rate is greater or equal than
			the maximal exit-rate $\alpha_0(x,t')=\sum_{j=1}^m \alpha_j(x,t')$
			over all states $x$ in $S_{t,R^*}$. In principal it is enough to
			find a function $\Lambda(t')$ with this property, for instance the
			function $\max_{x \in S_{t,R^*}} \sum_{j=1}^m \alpha_j(x,t')$, but
			this function may be hard to determine analytically and it is also
			not clear how to represent such a function practically in an
			implementation. Selecting a state $x^{\max}$ and defining $\Lambda(t')$
			to be the exit-rate of this state solves these problems.
			
			%  In Example~\ref{ex:2} we found that the maximal exit-rate
			% occurred in the same state throughout this time-interval, because the
			% only time-dependent transition class has the same rate in every state.
			% 
			% In general this is not the case and then we must determine, for each
			% time-point in $[t,t+\Delta)$ which state has the maximum exit-rate. It may
			% be possible to solve this problem analytically, but since the number of
			% states in $S_{t,R^*}$ may be large and this set is generated on-the-fly
			% we consider this to be infeasible.
			We now present two ways of implementing the function $\textrm{FindMaxState}$.

			%\setdefaultleftmargin{1.5em}{1em}{}{}{}{}
			\begin{itemize}
			\item[a)] For this approach we assume that all rate functions increase
				monotonically in the state variables. This is, for instance, always the
				case for models from chemical kinetics. We exploit that the change vectors are constant and define 
				for each dimension $k\in\{1,\ldots,n\}$
				$$w_k^{\max} :=\max_{j\in\{1,\ldots,m\}} w_{jk}$$
				where $w_{jk}$ is the $k$-th entry of the change vector $w_j$.
				For the set $S_{t,0}$ we compute, the maximum value for each dimension
				$k\in\{1,\ldots,n\}$
				$$ y_k^{\max} := \max_{y \in S_{t,0}} y_k. $$
			
				We now find the state $x^{\max}$ which is guaranteed to have
				a higher exit-rate than any state in $S_{t,R^*}$ for all time-points
				in the interval $[t,t+\Delta)$ as follows,
				$$
				x_k^{\max} := y_k^{\max} + R^* \cdot w_k^{\max}.
				$$
				It is obvious that the state variables $x_k^{\max}$ are upper bounds
				for the state variables appearing in $S_{t,R^*}$. Then, since all rates
				increase monotonically in the state variables, we have that the exit-rate
				of $x^{\max}=(x_1^{\max},\ldots,x_n^{\max})$ must be an upper-bound 
				for the exit-rates appearing in $S_{t,R^*}$ for all time-points.

			\item[b)] The first two moments of a Markov population model 
				can be accurately approximated using the method of moments 
				proposed by Engblom~\cite{engblom3}.  
				This approximation assumes that the expectations and the 
				(co-)variances change continuously and deterministically in time 
				and it is accurate for most models with rate functions that are at most
				quadratic in the state variables.
				We approximate the means $E_k(t') := E[X_k(t')]$ and the variances
				$\sigma_k^2(t') := \mathit{VAR}[X_k(t')]$ for all $k\in\{1,\ldots,n\}$. 
				For each $k$, we determine the time instant $\hat t\in[t,t+\Delta)$ at which 
				$E_k(\hat t)+ \ell \cdot \sigma_k(\hat t)$ is maximal
				for some fixed $\ell$.
				%  and define $\mu_k=E[X_k(\hat t)]$ and  
				% $\sigma^2_k=\mathit{VAR}[X_k(\hat t)]$. 
				We use this maximum to  determine the spread of the 
				distribution, i.e. we assume that the values of $X(t')$ will stay below 
				$x^{\max}_k:=E_k(\hat t)+ \ell \cdot \sigma_k(\hat t)$ with high probability.
				Note that a more detailed 
				approach is to consider the multivariate normal distribution 
				with mean $E[X(t')]$ and covariance matrix $\mathit{COV}[X(t')]$.
				But since the spread of a multivariate normal distribution 
				is difficult to derive in  higher dimensions, we simply
				consider each dimension independently. We now have $x^{\max} =
				(x^{\max}_1,\ldots,x^{\max}_n)$.
				 If during the analysis a state is found which exceeds $x^{\max}$ in
				one dimension then we repeat our computation 
				with a higher value for $\ell$. 
				To make this approach efficient,  $\ell$ has to 
				be chosen in an appropriate way. Our experimental results 
				indicate that for two-dimensional systems the choice $\ell=4$
				yields the best results.

			\end{itemize}

		\subsection{Complete algorithm}
			%Our complete algorithm now proceeds as follows (see Algorithm~\ref{alg:alg_pipeline}).
			Our complete algorithm now proceeds as follows (see Algorithm~\ref{alg:complete_algorithm}).
			Given an initial distribution
			$p^{(0)}$ with finite support $S_{0,0}$, a time-bound $t^{\max}$,
			thresholds $\delta$ and $\epsilon$, and
			a desired right truncation point $R^*$, we first set $t := 0$.
			Now we compute a time-step $\Delta$ and the state $x^{\max}$ using
			%Algorithm~\ref{alg:binarysearch} with inputs $R^*$, $t$, $t^{\max}$, and
			Algorithm~\ref{alg:binarysearch} with inputs $R^*$, $t$, $t^{\max}$, and
			%Algorithm~1 (we refer to it as ChooseTimeStep function) with inputs $R^*$, $t$, $t^{\max}$, and
			$\epsilon$. We then approximate the transient distribution $\hat p^{t+\Delta}$
			using an on-the-fly version of the bounding approach~\cite{buchholz}, where
			the state space is dynamically maintained and states with probability less
			than $\delta$ are discarded as described above. For the  rate function $\Lambda(t)$ we
			use the exit-rate of state $x^{\max}$. When computing DTMC probabilities,
			we use exact formulas for the first two terms~\cite{buchholz}
			of the sum in Eq.~\eqref{eq:approx} and lower bounds, given by Eq.~\eqref{eq:minbound},
			for the rest. This gives us the approximation $\hat p^{t+\Delta}$ with
			finite support $S_{t+\Delta,0}$. We now set $t := t+\Delta$ and repeat
			the above step with initial distribution $\hat p^{t}$ until we have $t = t^{\max}$.
			
	\section{Case Studies}\label{sec:expr}
		We implemented the approach outlined in Section~\ref{sec:onthefly}
		in C++ and ran experiments on a 2.4GHz
		Linux machine with 4 GB of RAM.
		We consider a Markov population model that describes a network 
		of chemical reactions.
		According to the theory of stochastic chemical kinetics~\cite{gillespie77},
		 the form of the rate function of a   reaction depends on 
		how many molecules of each chemical species are needed for 
		one instance of the reaction to occur. 
		The relationship to the volume has been discussed in detail 
		by Wolkenhauer et al.~\cite{wolkenhauer04}.
		If no reactants are needed\footnote{Typically, 
		reactions requiring no reactants  are used in the case of open systems where
		it is assumed that the reaction is always possible at a constant rate and 
		the reactant population is not explicitly modelled.}, that is,
		the reaction is of the form $\emptyset\to\ldots$ then 
		$\alpha_j(x,t)=k_j\cdot V(t)$ where $k_j$ is a positive constant
		and $V(t)$ is the volume of the compartment in which the reactions take 
		place. If one molecule is needed (case $S_i\to \ldots$) then 
		$\alpha_j(x,t)=k_j\cdot x_i$ where $x_i$ is the number of molecules 
		of type $S_i$. Thus, in this case, $\alpha_j(x,t)$ is independent of time.
		If two distinct molecules are needed (case $S_i+S_\ell\to\ldots)$
		then $\alpha_j(x,t)=\frac{k_j}{V(t)}\cdot x_i\cdot x_\ell$.

		All these theoretical considerations are based on the 
		assumption that the chemical reactions are \emph{elementary},
		that is, they are not a combination of several reactions.
		Our example may contain non-elementary reactions and 
		thus  a realistic biological model may contain different volume 
		dependencies. But since the focus of the paper is on the 
		numerical algorithm, we do not aim for an accurate 
		biological description here.

		\begin{table*}[h]
		\renewcommand{\arraystretch}{1.2}
			\begin{center}
				\begin{tabular}{|c|c|c|c|r|c|c|c|c|} 
					\hline
						Case study &
						\begin{tabular}{c} FindMaxState \\ implementation \end{tabular}		
						& $R^*$ 
						& \begin{tabular}{c} Total \\ error \end{tabular}
						%& \begin{tabular}{c} Execution \\ time \end{tabular}
						& Ex. time
						& $|S|$ 
						& $\min_{\%}$ 
						& $ \textrm{Poisson}_{\%}$ \\
					\hline
						\multirow{8}{*}{ \begin{tabular}{c} Gene \\ expression \end{tabular} }
		 				 & \multirow{4}{*}{method a)} &5 &$4.69 \cdot 10^{-4}$ & 14 min & \multirow{4}{*}{33962} 	& 95 & 5		\\
		 				 &   & 10 & $1.33 \cdot 10^{-2} $ & 10 min 	&  & 78 & 22  \\
		 				 &   & 15 & $2.24 \cdot 10^{-2} $ & 5 min	&  & 64 & 36  \\
						 & 	 & 20 & $9.92 \cdot 10^{-2}$ &  3 min 	&  & 41 & 59 \\
					\cline{2-8}
						& \multirow{4}{*}{method b)}	& 5	& $4.78 \cdot 10^{-4}$	&	27 min	&	\multirow{4}{*}{33130}& 95	& 5	\\
						&			& 10 	& $9.63 \cdot 10^{-3}$	&	14 min	&			&	77	& 23	\\
						&			& 15 	& $4.08  \cdot 10^{-2}$	&	10 min	&			&	58	& 42	\\
						&			& 20	& $7.73 \cdot 10^{-2}$	&	7 min	&			&	41	& 59	\\
			
					\hline
						\multirow{8}{*}{ \begin{tabular}{c} Exclusive \\ switch \end{tabular} }
						& \multirow{4}{*}{method a)}  	& 5  &$2.38 \cdot 10^{-6}$  &   21 min   & \multirow{4}{*}{1740}  & 80 & 20	\\
						&  	& 10 &$1.63 \cdot 10^{-5}$  & 29 min &     	& 75 & 25	\\
						&  	& 15 &$2.51 \cdot 10^{-5}$	& 68 min &    	& 47 & 53	\\
						&  	& 20 &$3.32 \cdot 10^{-5}$  & 2 h  	 &   	& 38 & 62	\\
					\cline{2-8}
						& \multirow{3}{*}{method b)}   	& 5 & $3.56 \cdot 10^{-6}$	& 17 h 		& \multirow{4}{*}{1752}  &89 	& 11 \\
						&				&10	& $1.55 \cdot 10^{-4}$ 	&  3 h 		& 	& 78	& 22 	\\
						&				&15	& $6.51 \cdot 10^{-4}$ 	&  1.5 h	& 	& 59	& 41 	\\
						&				&20	& $1.71 \cdot 10^{-3}$	& 1 h		& 	& 42	& 58	\\						
				    \hline
		 		\end{tabular}
			\end{center}
			\caption{Results of the analysis of case studies.\label{tab:switch}}
		\end{table*}

		%\subsection{Exclusive Switch}
		%Our second example is a 
		We conduct experiments with two reaction networks.
		The first one is a simple gene expression (described in Ex.~\ref{ex:gene}).
		The second one is a
		%The reaction network that we consider is a
		gene regulatory
		network, called the exclusive switch~\cite{exswitch}.
		It consists of two genes with a common promotor region.
		Each of the  two gene products $P_1$ and $P_2$ inhibits the 
		expression of the other product if a molecule is bound to the 
		promotor region. More precisely, 
		if the promotor region is free, molecules of  both types 
		$P_1$ and $P_2$ are produced.
		If a molecule of type $P_1$ is bound to the promotor 
		region, only molecules of    type $P_1$ are produced.
		If a molecule of type $P_2$ is bound to the promotor 
		region, only molecules of    type $P_2$ are produced.
		No other configuration of the promotor region exists.
		The probability distribution of the exclusive switch is bistable
		which means that after a certain amount of time, the 
		probability mass concentrates on two distinct regions 
		in the state space. 
		The system has five chemical species of which two 
		have an infinite range, namely $P_1$ and 
		$P_2$. We define the  transition classes 
		%$\tau_j=\left(G_j,u_j,\alpha_j\right)$, $j\in\{1,\ldots,10\}$ as follows.
		$\tau_j=\left(G_j,w_j,\alpha_j\right)$, $j\in\{1,\ldots,10\}$ as follows.
		 	\begin{itemize}
				 \item For $j\in\{1,2\}$ we describe production of  $P_j$
					by $G_j=\{x\in\mathbb N^5\mid x_3>0\}$, 
					$w_j=e_j$, 
					%$u_j(x)=x+e_j$, 
					and 
					$\alpha_j(x,t)=0.5 \cdot x_3$. Here, $x_3$ denotes the number 
					of unbound DNA molecules which is either zero or one and the vector  $e_j$
					is such that all its entries are zero except the $j$-th entry which is one.
				\item We describe degradation of  $P_j$
					by $G_{j+2}=\{x\in\mathbb N^5\mid x_j>0\}$, 
					$w_{j+2}=-e_j$,
					%$u_{j+2}(x)=x-e_j$,
					 and 
					$\alpha_{j+2}(x,t)=0.005 \cdot x_j$. Here, $x_j$ denotes the number 
					of $P_j$ molecules.
				\item We model the binding of $P_{j}$  to the promotor as
					$G_{j+4}=\{x\in\mathbb N^5\mid x_3>0,x_{j}>0\}$,
					$w_{j+4}=-e_j-e_3+e_{j+3}$, and 
					%$u_{j+4}(x)=x-e_j-e_3+e_{j+3}$, and 
					$\alpha_{j+4}(x,t)=(0.1-\frac{0.05}{3600}\cdot t) \cdot x_j\cdot x_3$
					for $t\le 3600$. Here, $x_{j+3}$ is one if a molecule of type $P_j$ 
					if bound to the promotor region and zero otherwise.
				\item For unbinding of $P_j$ we define
		 			$G_{j+6}=\{x\in\mathbb N^5\mid  x_{j+3}>0\}$,
		 			$w_{j+6}=e_j+e_3-e_{j+3}$, and 
					%$u_{j+6}(x)=x+e_j+e_3-e_{j+3}$, and 
					$\alpha_{j+6}(x,t)=0.005\cdot x_{j+3}$.
				 \item Finally, we have production of $P_j$ if a molecule
					of type $P_j$ is bound to the promotor, i.e.,
					$G_{j+8}=\{x\in\mathbb N^5\mid  x_{j+3}>0\}$,
					$w_{j+8}=e_j$, and 
					%$u_{j+8}(x)=x+e_j$, and 
					$\alpha_{j+8}(x,t)=0.5\cdot x_{j+3}$.
		 \end{itemize}
		%Note that only the rate functions $\alpha_6$ and $\alpha_7$,
		Note that only the rate functions $\alpha_5$ and $\alpha_6$,
		which denote the binding of a protein to the promotor region,
		are time-dependent. This is intuitively clear since 
		if the cell volume grows it becomes less likely that a protein 
		molecule is located close to the promotor region. 
		We started the system at time $t=0$ in state $(0,0,1,0,0)$ with probability 
		one and considered a time horizon of $t=3600$. 
		For the simple gene expression system (Example~\ref{ex:gene}) we started
		at time $t=0$ in state $(0,0)$ and considered the same time horizon.
		Table~\ref{tab:switch} contains the results of our experiments.
		The first column refers to the system under study and the second one
		shows the variation used to implement the method \textrm{FindMaxState} which we suggest in 
		Section~\ref{sec:maxrates}. 
		%The first column refers to the two variants for 
		%implementing the  method \textrm{FindMaxState} which we suggest in 
		%Section~\ref{sec:maxrates}. 
		%The second and third column lists the different 
		%values that we used for the threshold $\delta$ and the right truncation 
		%point $R^*$. 
		The third column lists the different values for right truncation
		point $R^*$. 
		We list the total error at time $\tmax$ in the fourth column
		(see Eq.~\eqref{eq:totalerror}). Program execution time 
		is given in the fifth column and the sixth column with heading
		$|S|$ contains the maximal size of the set $S_{t,R^*}$ that 
		we considered during the analysis. The next two columns describe
		the percentage of the total probability loss due to the
		bounding approach ($\min_{\%}$) and due to the truncation
		of the infinite sum in Eq.~\eqref{eq:unif_sum} ($\textrm{Poisson}_{\%}$).
		The two percentages in one row do not sum up to one
		since we store only states that have significant probability 
		(w.r.t threshold $\delta$), which is the third error source. 
		However, this lost portion is negligible for the two systems
		that we consider.
		For our implementation we kept the input $\epsilon=10^{-10}$ 
		%of Algorithm~\ref{alg:binarysearch} fixed.
		of Algorithm~\ref{alg:binarysearch} fixed.

		\subsection{Discussion}
			We now discuss the effect of the different input parameters on the
			performance of our algorithm and start with the implementation
			of the method to approximate $x^{\max}$. For both systems the method ''b''
			is less effective than method ''a'' (see Section~\ref{sec:maxrates}).
			Method ''b'' gives larger uniformization rates than method ''a'', which leads to
			slower execution times. Notice that the execution time grows
			when we use method ''a'' for the exclusive switch system when
			we choose the larger values for $R^*$. This is due to the fact
			that it always finds a state $x^{\max}$
			without taking expectations and covariances into consideration.
			This results in large over-approximations for such a bi-stable system. 
			The effect of the choice between methods ''a'' and ''b''
			on the accuracy is not completely clear, both methods provide
			the same order of the probability loss for the simple gene 
			expression system. For the second case study method ''a''
			provides tighter error bounds for larger values of $R^*$.
			
			In the Table~\ref{tab:switch} we show results obtained with
			$\delta=10^{-15}$. 
			%Naturally, choosing the lower threshold for
			%the significant probability mass (see Section~\ref{sec:onthefly}) 
			Naturally, choosing a lower threshold 
			results in larger execution times but one can gain a deeper exploration
			of the state space. This fact can also be used to obtain
			a coarse solution for certain system by setting $\delta=10^{-5}$,
			for instance.
			
			%As expected, decreasing the threshold
			%$\delta$ increases the accuracy, since less states are discarded on-the-fly.
			%However, this comes at a cost of using more memory, since more states
			%have to be represented, and the running time is also increased.
			
			%We also see that using method ``b'' to find the uniformization rate
			%is less effective than method ``a'' (see Section~\ref{sec:maxrates}). Method
			%``b'' chooses a larger uniformization rate than method ``a'', which leads to
			%slower execution times and increased memory usage. The effect of this
			%choice on the accuracy is not completely clear, although also here method
			%`a'' seems to be somewhat better. 
			
			The effect of the choice of $R^*$ is most interesting. Choosing a larger
			value for $R^*$ means that more summands on the right-hand side of
			Eq.~\eqref{eq:minbound} have to be approximated using the bounding approach.
			This decreases the accuracy of the algorithm since the larger time
			steps $\Delta$ are conducted and one obtain coarse approximation. 
			However it reduces the running time since $\tmax$
			can be covered using fewer iterations.
			Notice that the percentage of the
			probability loss due to truncation of the infinite sum in
			Eq.~\eqref{eq:unif_sum} grows when $R^*$ is chosen to be large.
			The reason is that we compute only first $3$ exact terms in the sum
			and remaining terms are approximations.	Thus the choice of $R^*$
			determines the compromise between running time and accuracy.
			
%			This should decrease the accuracy of the algorithm, but we see that for one
%			of the experiments this is not the case (method ``a'', $\delta=10^{-12}$).
%			This may be caused by the fact that increasing $R^*$ also increases
%			the time-steps $\Delta$ and the uniformization rate $\Lambda(t)$. By increasing
%			the time-steps we find that less steps have to be taken to reach the final
%			time-point $t^{\max}$ which decreases the probability lost by the truncation
%			of the uniformization sum. We also see that increasing $R^*$ increases
%			the memory and time needed for computation.

	\section{Conclusion}
		We have presented an algorithm for the numerical approximation of
		transient distributions for infinite time-inhomogeneous Markov
		population models with unbounded rates. Our algorithm provides a strict
		lower bound for this transient distribution. There is a trade-off between
		the tightness of the bound and the performance of the algorithm, both
		in terms of computation time and required memory.
		
		As future work, we will investigate the relationship between  the 
		parameters of our approach (truncation 
		point, the significance threshold $\delta$, the method by which we determine 
		the rate of the Poisson process), the accuracy  and the running time
		of the algorithm more closely. For this we will consider Markov population models
		with different structures and dynamics.

\balance

%  \bibliographystyle{plain}
% \bibliography{cav}

%\nocite{*}
\bibliographystyle{eptcs}
\bibliography{cav}

\end{document}